\documentclass[12pt]{amsart}
\usepackage{amsmath,amsthm,amssymb,amscd,latexsym,eufrak}

\newcommand{\svskip}{\vspace{3mm}}

\newcommand{\R}{{\Bbb R}}

\newcommand{\Z}{{\Bbb Z}}

\newcommand{\Ker}{{\rm Ker}\:}

\newtheorem{thm}{Theorem}[section]
\newtheorem{lem}[thm]{Lemma}

\newtheorem{problem}[thm]{{\sc Problem}}

\newtheorem{defn}[thm]{{\sc Definition}}

\begin{document}
\title[A generalization of Roberts' counterexample]{A note on "a generalization of Roberts' counterexample to the fourteenth problem of Hilbert by S.\ Kuroda"}
\author{Mikiya Tanaka}
\address{Graduate School of Science and Technology, Kwansei Gakuin University, 
2-1 Gakuen, Sanda, Hyogo 669-1337, Japan}
\email{mtanaka@kwansei.ac.jp}
\keywords{locally nilpotent derivation, the fourteenth problem of Hilbert}
\subjclass[2000]{Primary: 14R20; secondary: 13A50, 13N15}
\maketitle

\begin{abstract}
In \cite{Kuroda}, Kuroda generalized Roberts' counterexample \cite{Roberts} to the fourteenth problem of Hilbert.
The counterexample is given as the kernel of a locally nilpotent derivation on a polynomial ring.
We replace his construction of the invariant elements by a more straightforward construction and give a more precise form of invariant elements.
\end{abstract}

\section{Introduction}

Let $k$ be a field of characteristic zero and let $B$ be a $k$-algebra.
We denote by ${\rm LND}_k(B)$ the set of $k$-derivations of $B$.
In \cite{Kuroda}, Kuroda proved the following result.

\begin{thm}\label{thm1.1}
Let $B=k[x_1,\dots,x_n,y_1,\dots,y_n,y_{n+1}]$ be a polynomial $k$-algebra and define $\delta\in {\rm LND}_k(B)$ by $\delta(x_i)=0$ 
and $\delta(y_i)={x_i}^{2}$ for all $1\leq i\leq n$, and $\delta(y_{n+1})=x_1\cdots x_n$.
Suppose that $n\geq 4$.
Then $A:=\Ker \delta$ is not finitely generated over $k$.
\end{thm}

In order to prove this theorem, he made use of the following lemma.

\begin{lem}\label{lem1.2}
With the notations and assumptions in the above theorem, there exists a positive integer $\alpha$ such that the $k$-subalgebra $A$ contains elements of the form
$$
{x_1}^\alpha{y_{n+1}}^\ell+(\text{terms of lower degree in }y_{n+1})
$$
for each $\ell \geq 1$.
\end{lem}

In this paper, we prove that we can take $\alpha=1$.
Namely, we prove the following.

\begin{thm}\label{thm1}
With the notations and assumptions in the above theorem, the $k$-subalgebra $A$ contains elements of the form
$$
x_1{y_{n+1}}^\ell+(\text{terms of lower degree in }y_{n+1})
$$
for each $\ell \geq 1$.
\end{thm}

\section{Proof of Theorem \ref{thm1}}

In a subsequent proof, we use the following result.

\begin{lem}
Let $B=k[x_2,\dots,x_n,y_2,\dots,y_n]$ and define $\delta \in {\rm LND}_k(B)$ by $\delta(x_i)=0$ and $\delta(y_i)={x_i}^2$ for each $i$.
Then $\Ker \delta$ is a $k$-algebra generated by $x_2,\dots,x_n$ and ${x_i}^2y_j-{x_j}^2y_i\ (2\leq i,j\leq n, i\neq j)$.
\end{lem}

We can prove this lemma by the same argument in \cite[Theorem 1.2]{Kojima-Miyanishi}.

Now, for each monomial $m={x_1}^{a_1}\cdots{x_n}^{a_n}{y_1}^{b_1}\cdots{y_{n+1}}^{b_{n+1}}$, define
$$
\tau(m)=\left[\frac{a_2}{2}\right]+\cdots+\left[\frac{a_n}{2}\right]-(b_1+\cdots+b_n),
$$
where we write $[a]=\max\{n\in \Z\mid n\leq a\}$ for any $a\in \R$.
Let $f_{1,n+1}=x_1y_{n+1}-x_2\cdots x_ny_1$ and let $f_{i,j}={x_i}^2y_j-{x_j}^2y_i$ for each pair $(i,j)$ with $1\leq i,j\leq n$ and $i\neq j$.
It is easy to see that all of $f_{1,n+1}$, $x_i$ and $f_{i,j}$ belong to $A$.
Let $A'$ be the $k$-subalgebra generated by $f_{1,n+1}$, $x_i\ (1\leq \i \leq n)$ and $f_{i,j}\ (1\leq i,j \leq n, i\neq j)$.
Since $A$ is factorially closed in $B$, i.e., $a=b_1b_2$ with $b_1,b_2\in B$ implies $b_1,b_2\in A$, it suffices to show that there exists $f\in A'$ such that ${f_{1,n+1}}^\ell -f$ is of the form
\begin{align*}
&{x_1}^\ell {y_{n+1}}^\ell +{x_1}^{\ell-1}(\text{terms of lower degree in }y_{n+1})\\
&={x_1}^{\ell-1}\left(x_1{y_{n+1}}^{\ell}+(\text{terms of lower degree in }y_{n+1})\right).
\end{align*}
We have
\begin{align*}
{f_{1,n+1}}^\ell=&{x_1}^\ell{y_{n+1}}^\ell-\ell {x_1}^{\ell-1}{y_{n+1}}^{\ell-1}x_2\cdots x_n y_1\\
&+\binom{\ell}{2}{x_1}^{\ell-2}{y_{n+1}}^{\ell-2}{x_2}^2\cdots{x_n}^2{y_1}^2+\cdots+(-1)^\ell {x_2}^\ell \cdots {x_n}^\ell {y_1}^\ell
\end{align*}
and we construct $f\in A'$ which, when subtracted from ${f_{1,n+1}}^\ell$, cancels the terms in ${f_{1,n+1}}^\ell$ of degree $<\ell-1$ in $x_1$ and produces only the terms of degree $\geq \ell-1$ in $x_1$.
Namely, as the element ${f_{1,n+1}}^\ell-f$, we construct an element in $A'$ of the form
$$
{x_1}^\ell {y_{n+1}}^\ell + g_{l-1}{y_{n+1}}^{\ell-1}+g_{\ell-2}{y_{n+1}}^{\ell-2}+\cdots+g_0,
$$
where $g_i\in k[x_1,\dots,x_n,y_1,\dots,y_n]$ and ${x_1}^{\ell-1}$ divides every $g_i$.

By the descending induction on $r$, we suppose that we obtain an element in $A'$ of the form
$$
G_{r}={x_1}^\ell {y_{n+1}}^\ell + g_{l-1}{y_{n+1}}^{\ell-1}+g_{\ell-2}{y_{n+1}}^{\ell-2}+\cdots+g_r{y_{n+1}}^r+\cdots+g_0
$$
with $g_i\in k[x_1,\dots,x_n,y_1,\dots,y_n]$ and $g_{\ell-1},\dots,g_{r+1}$ divisible by ${x_1}^{\ell-1}$.
We show that $G_r$ is modified by an element of $A'$ so that a new $g_r$ is divisible by ${x_1}^{\ell-1}$ without changing the terms $g_{\ell-1},\dots,g_{r+1}$.
Furthermore, we suppose the following conditions are satisfied.
\begin{enumerate}
\item[(1)]For $0\leq i \leq \ell-1$, if we write $g_i=\sum_j {x_1}^j {y_1}^{q_{i,j}} h_{i,j}$ with $h_{i,j}\in k[x_2,\dots,x_n,y_2,\dots,y_n]$, then $i+j+2q_{i,j}=2\ell$.
\item[(2)]We have $h_{i,0}=\cdots=h_{i,i-1}=0$ for $0\leq i\leq \ell-1$, i.e., for each ${x_1}^j {y_1}^{q_{i,j}} h_{i,j}$ appearing in $g_i$, we have $j\geq i$.
\item[(3)]For each monomial $m={x_1}^j{x_2}^{a_2}\cdots{x_n}^{a_n}{y_1}^{b_1}\cdots{y_n}^{b_n}{y_{n+1}}^i$\\
in ${y_{n+1}}^i{x_1}^j {y_1}^{q_{i,j}} h_{i,j}$, we have
\begin{enumerate}
\item[(i)]$2\tau(m)\geq \ell-j-3$ and $a_2,\dots,a_n$ are all odd integers if $j \equiv \ell-1 \pmod{2}$,
\item[(ii)]$2\tau(m)\geq \ell-j$ and $a_2,\dots,a_n$ are all even integers if $j\equiv \ell \pmod{2}$.
\end{enumerate}
\end{enumerate}
In order to improve the term $g_r$ in such a way that $h_{r,0}=\cdots=h_{r,\ell-2}=0$, we suppose by a double induction that $h_{r,0}=\cdots=h_{r,p-1}=0$ and $h_{r,p}\neq 0$ with $r\leq p \leq \ell-2$.
With this hypothesis taken into account, we denote the polynomial $G_r$ by $G_{r,p}$.
The beginning polynomial for induction is $G_{\ell-2,\ell-2}={f_{1,n+1}}^\ell$, for which $g_i=(-1)^i \binom{\ell}{i}{x_1}^i(x_2\cdots x_n y_1)^{\ell-i}$, $h_{i,i}=(-1)^i \binom{\ell}{i}(x_2\cdots x_n)^{\ell-i}$, $h_{i,j}=0\ (i\neq j)$ and $q_{i,i}=\ell-i$ for $0\leq i,j\leq \ell-1$.
One can check easily that the above conditions are satisfied for $G_{\ell-2,\ell-2}={f_{1,n+1}}^\ell$

We explain the process of improving $g_r$.
Since $g_{r+1}$ is divisible by ${x_1}^{\ell-1}$ and 
\begin{align*}
&\delta({x_1}^p{y_1}^{q_{r,p}}h_{r,p}{y_{n+1}}^r)={x_1}^p{y_1}^{q_{r,p}}{y_{n+1}}^r\delta(h_{r,p})\\
&+(\text{terms of degree }>p\text{ in }x_1)+(\text{terms of degree }<r\text{ in }y_{n+1}),
\end{align*}
we have
\begin{align*}
&0=\delta(G_{r,p})={x_1}^p{y_1}^{q_{r,p}}{y_{n+1}}^r\delta(h_{r,p})\\
&+(\text{terms of degree }>p\text{ in }x_1)+(\text{terms of degree }\neq r\text{ in }y_{n+1})
\end{align*}
and hence $\delta(h_{r,p})=0$.
Lemma 2.1 implies that $h_{r,p}$ is a sum of polynomials of the form
$$
c{x_2}^{d_2}\cdots{x_n}^{d_n}\prod_{i,j\in \{2,\dots,n\}}{f_{i,j}}^{t_{i,j}}
$$
with $c\in k$ and non-negative integers $d_i$, $t_{i,j}$.
Note that all of $d_2,\dots,d_n$ are odd integers (resp.\ even integers) if $p\equiv \ell-1 \pmod{2}$ (resp.\ if $p\equiv \ell\pmod{2}$).
In fact, since the contributions of the $f_{i,j}$ to the exponent $d_2,\dots,d_n$ are even, the remark follows from the conditions (i) and (ii) of (3).
Now we choose any one of the above polynomials and let $H=\prod_{i,j\in \{2,\dots,n\}}{f_{i,j}}^{t_{i,j}}$.
Then, for each monomial $m$ in ${y_{n+1}}^r{x_1}^p{y_1}^{q_{r,p}}{x_2}^{d_2}\cdots{x_n}^{d_n}H$, we have in view of (i) and (ii) of (3),
\begin{align*}
2\tau({y_1}^{q_{r,p}}{x_2}^{d_2}\cdots{x_n}^{d_n})&=2\tau(m)\\
&\geq 
\begin{cases}
\ell-p-3 & \text{if } p\equiv \ell-1 \pmod{2}\\
\ell-p & \text{if } p\equiv \ell \pmod{2}
\end{cases},
\end{align*}
where multiplying ${y_1}^{q_{r,p}}{x_2}^{d_2}\cdots{x_n}^{d_n}$ by any ${x_i}^2y_j\ (i\neq 1)$, $y_{n+1}$ or $x_1$ does not change the value of $\tau$.
Note that $r\leq p\leq \ell-2$ and that if $p\equiv \ell-1 \pmod{2}$, then $p\leq \ell-3$ and hence $\ell-p-3\geq 0$.
Thus we have $\tau({y_1}^{q_{r,p}}{x_2}^{d_2}\cdots{x_n}^{d_n})\geq 0$ and there exists an element $F\in A'$ of the form
\begin{align*}
F&=c{x_1}^{p-r}{f_{1,n+1}}^r{f_{2,1}}^{q_2}\cdots{f_{n,1}}^{q_n}{x_2}^{d_2-2q_2}\cdots{x_n}^{d_n-2q_n}\\
&=c{x_1}^p{y_1}^{q_{r,p}}{y_{n+1}}^r{x_2}^{d_2}\cdots{x_n}^{d_n}\\
&+(\text{terms of degree }>p\text{ in }x_1)+(\text{terms of degree }<r\text{ in }y_{n+1}),
\end{align*}
where $q_2+\cdots+q_n=q_{r,p}$.
We can prove that $G_{r,p}-FH$ satisfies the same conditions as $G_{r,p}$ does except for the condition $h_{r,p}\neq 0$ but the number of nonzero terms in $h_{r,p}$ gets smaller.
We prove this below.
By repeating this process finitely many times, we obtain a new $G_r$ satisfying the condition $h_{r,p}=0$.
Further, continuing this process finitely many times, we obtain a modified $G_r$ satisfying the condition $h_{r,0}=\cdots=h_{r,\ell-2}=0$, i.e., $g_r$ is divisible by ${x_1}^{\ell-1}$.
Hence by induction on $r$, we completes a proof.

Now we show that $G_{r,p}-FH$ satisfies the same conditions as $G_{r,p}$ does but the number of nonzero monomial terms in $h_{r,p}$ becomes less.
We have only to show that each monomial in $F$ satisfies the conditions (1)-(3) since none of $y_{n+1}$, $x_1$ and $y_1$ appears in $H$ and multiplication of any monomial in $H$ to a monomial does not chage the value of $\tau$.
Each nonzero monomial $m_F$ in $F$ is of the form
\begin{align*}
{x_1}^{p-r}(x_1y_{n+1})^{r_1}(x_2\cdots x_n y_1)^{r_2}{x_2}^{d_2-2q_2}\cdots{x_n}^{d_n-2q_n}
\prod_{i=2}^n ({x_i}^2 y_1)^{\alpha_i}({x_1}^2y_i)^{\beta_i}
\end{align*}
with $r_1+r_2=r$ and $\alpha_i+\beta_i=q_i$ for $i=2,\dots,n$.
We choose one $m_F$ and let $w$, $z_1$, and $z_{n+1}$ be the exponents of $x_1$, $y_1$, and $y_{n+1}$ in $m_F$ respectively.
Then we have
\begin{align*}
w&=p-r+r_1+2\beta_2+\cdots+2\beta_n=p-r_2+2\beta_2+\cdots+2\beta_n,\\
z_1&=r_2+\alpha_2+\cdots+\alpha_n \text{\quad and\quad}z_{n+1}=r_1.
\end{align*}
First we prove $m_F$ satisfies the conditions (1) and (2).
Indeed, we have
\begin{align*}
&z_{n+1}+w+2z_1=p+r_1+r_2+2(\alpha_2+\beta_2)+\cdots+2(\alpha_n+\beta_n)\\
&=p+r+2q_2+\cdots+2q_n=p+r+2q_{r,p}=2\ell
\end{align*}
and
\begin{align*}
&w-z_{n+1}=p-(r_1+r_2)+2(\beta_2+\cdots+\beta_n)\\
&=p-r+2(\beta_2+\cdots+\beta_n)
\geq p-r\geq 0.
\end{align*}

In order to prove that $m_F$ satisfies the condition (3), we consider four cases
\begin{enumerate}
\item[(a)]$p\equiv \ell-1 \pmod{2}$ and $r_2=2u+1$
\item[(b)]$p\equiv \ell-1 \pmod{2}$ and $r_2=2u$
\item[(c)]$p\equiv \ell \pmod{2}$ and $r_2=2u+1$
\item[(d)]$p\equiv \ell \pmod{2}$ and $r_2=2u$,
\end{enumerate}
where $u$ is an integer.
We only consider the case (a).
The remaining cases can be treated in a similar fashion.
Then we have
$$w\equiv \ell-1-2u-1+2\beta_2+\cdots+2\beta_n\equiv \ell \pmod{2}.$$
The exponent of each $x_i\ (i\neq 1)$ in $m_F$ is equal to $2u+1+d_i-2q_i+2\alpha_i$.
Since each $d_i$ is an odd integer by the condition (i) of (3), it is an even integer.
In addition, we have
\begin{align*}
2\tau(m_F)=&2((n-1)u+(n-1)+\tau({x_2}^{d_2-2q_2}\cdots{x_n}^{d_n-2q_n})\\
&-(2u+1)-(\beta_2+\cdots+\beta_n))\\
=&2u(n-3)+2(n-2)+2\tau({x_2}^{d_2-2q_2}\cdots{x_n}^{d_n-2q_n})\\
&-2(\beta_2+\cdots+\beta_n)\\
=&(r_2-1)(n-3)+2(n-2)+2\tau({y_1}^{q_{r,p}}{x_2}^{d_2}\cdots{x_n}^{d_n})\\
&-2(\beta_2+\cdots+\beta_n)\\
\geq &r_2-1+2\cdot 2+\ell-p-3-2(\beta_2+\cdots+\beta_n)\\
=&\ell-(p-r_2+2\beta_2+\cdots+2\beta_n)=\ell-w,
\end{align*}
where the term $(n-1)$ in the first equality is due to the condition that all the $d_i$ and $r_2$ are odd integers and we use the condition $n\geq 4$ to show the inequality.
Thus the condition (3) holds for $m_F$.
This induction completes a proof of Theorem \ref{thm1}.

\section{Application to module derivations}

In this section, we give application of Theorem \ref{thm1} to locally nilpotent module derivations.
First, we recall the following definition (see \cite{Tanaka}).

\begin{defn}{\em
Let $\delta\in {\rm LND}_k(B)$ and let $M$ be a $B$-module with a $k$-linear endomorphism $\delta_M:M\to M$.
A pair $(M,\delta_M)$ is called a $(B,\delta)$-module (a $\delta$-module, for short) if the following two conditions are satisfied.
\begin{enumerate}
\item[(1)]For any $b\in B$ and $m\in M$, $\delta_M(bm)=\delta(b)m+b\delta_M(m)$.
\item[(2)]For each $m\in M$, there exists a positive integer $N$ such that ${\delta_M}^n(m)=0$ if $n\geq N$.
\end{enumerate}
Let $A=\Ker \delta$.
Then $\delta_M$ is an $A$-module endomorphism.
Whenever we consider $\delta$-modules, the derivation $\delta$ on $B$ is fixed once for all.
We call $\delta_M$ a {\em module derivation} (resp.\ {\em locally nilpotent module derivation}) on $M$ if it satisfies the condition (1) (resp.\ both conditions (1) and (2)).
}\end{defn}

If there is no fear of confusion, we simply say that $M$ is a $\delta$-module instead of saying that $(M,\delta_M)$ is a $\delta$-module.
If $M$ is a $\delta$-module, then $M_0:=\Ker \delta_M=\{m\in M \mid \delta_M(m)=0\}$ is an $A$-module.
We retain below the notations $A,M_0$ for this specific purposes.
For the basic properties of $\delta$-modules, we refer the readers to \cite{Tanaka}.

We consider the following problem.

\begin{problem}\label{problem}
Let $B$ be an affine $k$-domain with a locally nilpotent derivation $\delta$ and let $M$ be a finitely generated $B$-module with $\delta$-module structure.
Is $M_0$ a finitely generated $A$-module?
\end{problem}

We have positive answers to Problem \ref{problem} if one of the following conditions is satisfied (see \cite{Tanaka2}).
\begin{enumerate}
\item[(i)]$M$ is torsion-free as a $B$-module and $A$ is a noetherian domain.
\item[(ii)]$M$ is torsion-free as a $B$-module and $\dim B\leq 3$.
\item[(iii)]$M_0$ is a free $A$-module.
\item[(iv)]The $B$-module $BM_0$ generated by $M_0$ is a free $B$-module with a basis $\{e_1,\dots,e_n\}$ such that $e_i\in M_0$.
\item[(v)]$B=A[y]$ is a polynomial ring over a noetherian domain $A$, $a:=\delta(y)$ is a nonzero element of $A$ and $a$ has no torsion in $M$.
\end{enumerate}
In \cite{Tanaka}, there is an easy counterexample to Problem \ref{problem} in the case where $M$ has torsion as a $B$-module.
In addition, there are counterexamples in the free case by making use of the counterexamples to the fourteenth problem of Hilbert given by Roberts \cite{Roberts}, Kojima-Miyanishi \cite{Kojima-Miyanishi}, Freudenburg \cite{Freudenburg} and Daigle-Freudenburg \cite{Daigle-Freudenburg}.
In such examples, we take $B$ to be a polynomial ring and $M$ to be the differential module $\Omega_{B/k}$.
We can give $\Omega_{B/k}$ a natural module derivation as follows.

\begin{lem}
Let $B$ be a $C$-algebra and let $\delta$ be a locally nilpotent $C$-derivation of $B$.
Then the differential module $M:=\Omega_{B/k}$ is a $\delta$-module if we define $\delta_M$ by $\delta_M(db)=d\delta(b)$ for $b\in B$.
\end{lem}

We can prove this lemma easily (see \cite{Tanaka2}).
Theorem \ref{thm1} gives a new counterexample to Problem \ref{problem}.
Namely, we have the following assertions.

\begin{thm}\label{thm}
With the notations and assumptions in Theorem \ref{thm1.1}, let $M=\Omega_{B/k}$ be the differential module with natural $\delta$-module structure.
Namely, $M$ is a free $B$-module
$$
M=\bigoplus_{i=1}^n Bdx_i\oplus \bigoplus_{i=1}^n Bdy_i
$$
with a free basis $\{dx_1,\dots,dx_n,dy_1,\dots,dy_{n+1}\}$ and a module derivation defined by
\begin{align*}
&\delta(dx_i)=0,\quad \delta(dy_i)=2x_idx_i\ (1\leq i\leq n)\quad \text{and}\\
&\delta(dy_{n+1})=\sum_{i=1}^n x_1\cdots \stackrel{\vee}{x_i} \cdots x_n dx_i.
\end{align*}
Then $M_0$ is not a finitely generated $A$-module.
\end{thm}

We can prove this in a fashion similar to \cite[Theorem6.2]{Tanaka2} by making use of Theorem \ref{thm1}.
The fact that we can take $\alpha=1$ in Lemma \ref{lem1.2} plays an important role in the proof of Theorem \ref{thm}.

\end{document}